\documentclass[12pt]{amsart}
\usepackage{amsmath, amssymb, latexsym, amsthm, mathrsfs, color, mathtools, latexsym}
\usepackage[hyphens]{url}
\usepackage[normalem]{ulem}
\usepackage[T1]{fontenc} %for the accent in Kakol's name
\usepackage[margin=1.2in]{geometry}

\newcommand{\nc}{\newcommand}

\nc{\thusfar}{\par\bigskip\centerline{\my{--- Edited thus far ---}}\par\bigskip}
\nc{\lei}{\le^\oo}
\nc{\card}[1]{\left|#1\right|}
\nc{\medcard}[1]{\biggl|\,#1\,\biggr|}
\nc{\smallcard}[1]{|\,#1\,|}
\nc{\bds}{bidirectional $\roth$-scale}
\nc{\bbT}{\mathbb{T}}
\nc{\bbN}{\mathbb{N}}%{\w}
\nc{\beq}{\begin{eqnarray*}}\nc{\eeq}{\end{eqnarray*}}
\nc{\mbq}{\mb{?}}
\nc{\mb}[1]{{\mbox{\textbf{#1}}}}
\nc{\nop}{$\times$}
\nc{\fbn}{\!\!\fbox{\!\nop\!}\!\!}
\nc{\yup}{\checkmark}
\nc{\forces}{\Vdash}
\nc{\name}[1]{\dot{#1}}
\nc{\tf}{\my{FINISHED THUS FAR}}
\nc{\FU}{Fr\'echet--Urysohn}
\nc{\gs}{$\gamma$~space}
\nc{\Ga}{\Gamma}\nc{\Om}{\Omega}
\nc{\smallbinom}[2]{\begin{psmallmatrix} #1\\ #2 \end{psmallmatrix}}
\nc{\bgamma}{\smallbinom{\Om}{\Ga}}
\newcommand{\two}{\{0,1\}}
\nc{\productive}[2]{\bigl(#1,\allowbreak #2\bigr)^\x}
\nc{\prdct}[1]{\bigl(#1\bigr)^\x}
\nc{\Sel}{\mathsf{S}}
\nc{\sset}[2]{\{\,#1 : #2\,\}}
\nc{\smb}[1]{{\!\!\mb{#1}\!\!}}
\nc{\medset}[2]{{\biggl\{\,#1 : #2\,\biggr\}}}
\nc{\smallmedset}[2]{{\bigl\{\,#1 : #2\,\bigr\}}}
\nc{\set}[2]{{\left\{\,#1 : #2\,\right\}}}
\nc{\seq}[2]{{\la\, #1 : #2\,\ra}}
\nc{\eseq}[1]{#1_1, \allowbreak #1_2, \allowbreak\dotsc} %explicit sequence
\nc{\cube}{(\Cantor)^\bbN}
\nc{\Match}{\op{Match}}
\nc{\concat}[1]{\hat{\phantom{a}}\langle #1\rangle}
\nc{\poset}{\mathbb{P}}
\nc{\fn}[1]{{\op{Fn}(#1\times\w,2)}}
\nc{\linadd}{\op{linadd}}
\nc{\nonprod}{\non^\x}
\nc{\alephes}{{\aleph_0}}
\nc{\my}[1]{\marginpar{\textcolor{red}{***}}\textcolor{red}{#1}}
\nc{\later}[1]{{\color{green} #1}}
\nc{\BTs}[1]{{\color{green} #1 (BT)}}
\nc{\Cp}{\op{C}_\mathrm{p}}
\nc{\Bp}{\op{B}_p}
\nc{\Pa}[8]{\bibitem{#1} {#2}, \emph{#3}, {#4} \textbf{#5} ({#6}), {#7}--{#8}.}
\nc{\tPa}[5]{\bibitem{#1} {#2}, \emph{#3}, {#4}, to appear.}
\nc{\sPa}[4]{\bibitem{#1} {#2}, \emph{#3}, {#4}, submitted.}
\nc{\Bc}[9]{\bibitem{#1} {#2}, \emph{#3}, in: \textbf{#4} (#5), #6 #7, #8--#9.}
\nc{\fD}{\mathfrak{D}}
\nc{\fX}{\mathfrak{X}}
\nc{\Onbd}{\Op_{\mathrm{nbd}}} %{\Op_{\mathsf{nbd}}}
\nc{\Omnb}{\Om_{\mathrm{nbd}}} %{\Om_{\mathsf{nbd}}}
\nc{\od}{\mathfrak{od}}
\nc{\Setting}[7]{\xymatrix@R=4pt@C=7pt{#1\ar@{-}[r]&#2\ar@{-}[r]&#3\\&#4\ar@{-}[u]\\
#5\ar@{-}[uu]\ar@{-}[r] & #6\ar@{-}[u]\ar@{-}[r] & #7\ar@{-}[uu]}}
\nc{\mx}[1]{\begin{matrix}#1\end{matrix}}
\nc{\plim}{p\txt{-}\lim}
\nc{\Bgp}{{\Z^\bbN}}
\nc{\Cgp}{{{\Z_2}^\bbN}}
\nc{\Cite}[1]{\textbf{[#1]}}
\nc{\Next}[1]{{#1^+}}
\nc{\cFin}{\mathrm{cF}}
\nc{\scsp}{\text{-scale space}}
\nc{\cfn}{\text{cofinal}\ }
\nc{\Con}{\text{Concentrated}}
\nc{\Lind}{\text{Lindel\"of}\,}
\nc{\con}{\text{-Concentrated}}
\nc{\lind}{\text{-Lindel\"of}\,}
\nc{\ctbl}{\text{countably}\,}
\nc{\Men}{\text{Menger}}
\nc{\men}{\text{-Menger}}
\nc{\Hur}{\text{Hurewicz}}
\nc{\intvl}[2]{{[#1(#2),\allowbreak #1(#2\!+\!1))}}
\nc{\Dfin}{\mathfrak{D}_\mathrm{fin}}
\nc{\grbl}{{\mbox{\textit{\tiny gp}}}}
\nc{\bbP}{\mathbb{P}}
\nc{\BOfat}{\B_{\Om_{\mathrm{fat}}}}%\B_{\mathrm{fat}}}
\nc{\Bgood}{\B_{\mathrm{good}}}
\nc{\compactN}{\cl{\mathbb{N}}}
\nc{\blocks}[2]{\op{cl}_{#2}(#1)}
\nc{\blocksplus}[2]{\op{cl}^+_{#2}(#1)}
\nc{\arx}[1]{\texttt{http://arxiv.org/math/#1}}
\nc{\bq}{\begin{quote}}
\nc{\eq}{\end{quote}}
\nc{\cl}[1]{\overline{#1}}
\nc{\CH}{the Continuum Hypothesis}
\nc{\MA}{Martin's Axiom}
\nc{\Bfat}{\B_\mathrm{fat}}
\nc{\inv}{^{-1}}
\nc{\Cantor}{{\two^\bbN}}%{{2^\w}}
\nc{\bP}{\mathbf{P}}
\nc{\bof}{\op{\fb}}
\nc{\dof}{\op{\fd}}
\nc{\bofF}{\bof(\cF)}
\nc{\sr}[3]{\underset{\mbox{#3}}{\mbox{#1}}}
\nc{\gp}{\binom{\Om}{\Ga}}
\nc{\gpsmall}{\mbox{$\gp$}}
\nc{\gig}{\gimel}%{\gimel\Ga}
\nc{\gns}{\sone(\Om,\gig)}
\nc{\nsr}[2]{#1}
\nc{\Srg}{{\mathbb{S}}}
\nc{\Srgs}{{\mathbb{S}^*}}
\nc{\NN}{{\bbN^{\bbN}}}
\nc{\ZN}{{\Z^{\bbN}}}
\nc{\NNup}{{\bbN^{\uparrow\bbN}}}
\nc{\Pof}{\op{P}}
\nc{\PN}{{\Pof(\bbN)}}
\nc{\roth}{{[\bbN]^{\mbox{\tiny $\infty$}}}} %{{[\w]^{\w}}}
\nc{\Fin}{[\bbN]^{\text{$<\!\!\infty$}}} %{{[\w]^{<\w}}}%{{[\bbN]^{<\aleph_0}}}
\nc{\ici}{[\bbN]^{ \infty, \infty}}%{{[\bbN]^{(\aleph_0,\aleph_0)}}}
\nc{\Inc}{{\compactN^{\uparrow\bbN}}}
\nc{\powInc}[1]{{\big(\Inc\big)^{#1}}}
\nc{\powFin}[1]{{\big(\Fin\big)^{#1}}}
\nc{\powPN}[1]{{\big(\PN\big)^{#1}}}
\nc{\NcompactN}{{\compactN^\bbN}}
\nc{\Uarrow}{\smash{\big\uparrow}}
\nc{\LE}{\preccurlyeq}
\nc{\GE}{\succcurlyeq}
\nc{\op}{\operatorname}
\nc{\im}{\op{im}}
\nc{\Span}{\op{span}}
\nc{\maxfin}{\op{maxfin}}
\nc{\ran}{\op{range}}
\nc{\iso}{\cong}
\nc{\Madd}{{\M}^\star}
\nc{\cI}{\mathcal{I}}
\nc{\cJ}{\mathcal{J}}
\nc{\scrA}{\mathscr{A}}
\nc{\scrB}{\mathscr{B}}
\nc{\scrC}{\mathscr{C}}
\nc{\scrD}{\mathscr{D}}
\nc{\scrF}{\mathscr{F}}
\nc{\scrK}{\mathscr{K}}
\nc{\A}{\forall}
\nc{\B}{\mathrm{B}}
\nc{\cB}{\mathcal{B}}
\nc{\bB}{\mathbf{B}}
\nc{\BS}{\mathbf{B}(\mathcal{S})}
\nc{\BF}{\mathbf{B}(\mathcal{F})}
\nc{\BU}{\mathbf{B}(\mathcal{U})}
\nc{\cSp}{\mathcal{S}^+}
\nc{\cFp}{\mathcal{F}^+}
\nc{\cUp}{\mathcal{U}^+}

\nc{\BG}{\B_\Ga}
\nc{\BL}{\B_\Lambda}
\nc{\BT}{\B_\Tau}
\nc{\BTstar}{\B_{\Tau^*}}
\nc{\BO}{\B_\Om}
\nc{\DO}{\cD_\Om}
\nc{\KO}{\cK_\Om}
\nc{\CG}{C_\Ga}
\nc{\CL}{C_\Lambda}
\nc{\CT}{C_\Tau}
\nc{\CTstar}{C_{\Tau^*}}
\nc{\CO}{C_\Om}
\nc{\COgp}{C_{\Om^{\grbl}}}
\nc{\CLgp}{C_{\Lambda^{\grbl}}}
\nc{\BOgp}{\B_{\Om}^{\grbl}}
\nc{\BLgp}{\B_{\Lambda^{\grbl}}}
\nc{\sfC}{\mathsf{C}}
\nc{\sfD}{\mathsf{D}}
\nc{\bD}{\mathbf{D}}
\nc{\Tau}{\mathrm{T}}
\nc{\cA}{\mathcal{A}}
\nc{\cK}{\mathcal{K}}
\nc{\cD}{\mathcal{D}}
\nc{\cF}{\mathcal{F}}
\nc{\cS}{\mathcal{S}}
\nc{\cT}{\mathcal{T}}
\nc{\cG}{\mathcal{G}}
\nc{\cY}{\mathcal{Y}}
\nc{\J}{\mathcal{J}}
\nc{\cL}{\mathcal{L}}
\nc{\cM}{\mathcal{M}}
\nc{\cN}{\mathcal{N}}
\nc{\cH}{\mathcal{H}}
\nc{\cO}{\mathcal{O}}
\nc{\Op}{\mathrm{O}}
\nc{\rmA}{\mathrm{A}}
\nc{\rmF}{\mathrm{F}}
\nc{\rmB}{\mathrm{B}}
\nc{\rmD}{\mathrm{D}}
\nc{\rmP}{\mathrm{P}}
\nc{\cC}{\mathcal{C}}
\nc{\cP}{\mathcal{P}}
\nc{\bbQ}{\mathbb{Q}}
\nc{\bbR}{\mathbb{R}}
\nc{\cU}{\mathcal{U}}
\nc{\Un}{\bigcup}
\nc{\cV}{\mathcal{V}}
\nc{\cW}{\mathcal{W}}
\nc{\Z}{{\mathbb Z}}
\nc{\Impl}{\Rightarrow}
\long\def\forget#1\forgotten{\marginpar{\textcolor{green}{Forgetting...}}}
\nc{\ft}{\mathfrak{t}}
\nc{\fb}{\mathfrak{b}}
\nc{\fc}{\mathfrak{c}}
\nc{\fd}{\mathfrak{d}}
\nc{\fg}{\mathfrak{g}}
\nc{\oo}{\infty}
\nc{\fr}{\mathfrak{r}}
\nc{\fk}{\mathfrak{k}}
\nc{\bidi}{\mathfrak{bidi}}
\nc{\fu}{\mathfrak{u}}
\nc{\fh}{\mathfrak{h}}
\nc{\fp}{\mathfrak{p}}
\nc{\fj}{\mathfrak{j}}
\nc{\fs}{\mathfrak{s}}
\nc{\w}{\omega}
\nc{\x}{\times}
\nc{\Iff}{\Leftrightarrow}

\nc{\nin}{\notin}
\nc{\cat}{\hat{\ }}
\nc{\sub}{\subseteq}
\nc{\spst}{\supseteq}
\nc{\sm}{\setminus}
\nc{\as}{\subseteq^*}%{\let\proclaim\relax}
\nc{\les}{\le^*}
\nc{\leinf}{\le^{\infty}}
\nc{\leS}{\le_S}
\nc{\leF}{\le_{\mathcal{F}}}
\nc{\leU}{\le_{\mathcal{U}}}
\nc{\rest}{\restriction}

\nc{\la}{\langle}
\nc{\ra}{\rangle}
\nc{\E}{\exists}
\nc{\dom}{\op{dom}}
\nc{\cov}{\op{cov}}
\nc{\add}{\op{add}}
\nc{\cof}{\op{cof}}
\nc{\cf}{\op{cf}}
\nc{\non}{\op{non}}
\nc{\unif}{\op{non}}
\nc{\COV}{\op{COV}}
\nc{\ADD}{\op{ADD}}
\nc{\COF}{\op{COF}}
\nc{\NON}{\op{NON}}
\nc{\impl}{\to}
\nc{\Lp}{\mathcal{L_\p}}
\nc{\Wlog}{without loss of generality}
\newtheorem{thm}{Theorem}[section]
\nc{\bthm}{\begin{thm}} \nc{\ethm}{\end{thm}}
\newtheorem{prop}[thm]{Proposition}
\nc{\bprp}{\begin{prop}} \nc{\eprp}{\end{prop}}
\newtheorem{fact}[thm]{Fact}
\nc{\bfct}{\begin{fact}} \nc{\efct}{\end{fact}}
\newtheorem{prob}[thm]{Problem}
\nc{\bprb}{\begin{prob}} \nc{\eprb}{\end{prob}}
\newtheorem{lem}[thm]{Lemma}
\nc{\blem}{\begin{lem}} \nc{\elem}{\end{lem}}
\newtheorem{app}[thm]{Application}
\nc{\bapp}{\begin{app}} \nc{\eapp}{\end{app}}
\newtheorem{cor}[thm]{Corollary}
\nc{\bcor}{\begin{cor}} \nc{\ecor}{\end{cor}}
\newtheorem{conj}[thm]{Conjecture}
\nc{\bcnj}{\begin{conj}} \nc{\ecnj}{\end{conj}}
\theoremstyle{definition}
\newtheorem{defn}[thm]{Definition}
\nc{\bdfn}{\begin{defn}} \nc{\edfn}{\end{defn}}
\newtheorem{obs}[thm]{Observation}
\nc{\bobs}{\begin{obs}} \nc{\eobs}{\end{obs}}
\theoremstyle{remark}
\newtheorem{rem}[thm]{Remark}
\nc{\brem}{\begin{rem}} \nc{\erem}{\end{rem}}
\newtheorem{claim}[thm]{Claim}
\nc{\bclm}{\begin{claim}} \nc{\eclm}{\end{claim}}
\newtheorem{cnv}[thm]{Convention}
\nc{\bcnv}{\begin{cnv}} \nc{\ecnv}{\end{cnv}}
\newtheorem{exam}[thm]{Example}
\nc{\bexm}{\begin{exam}} \nc{\eexm}{\end{exam}}
\nc{\bpf}{\begin{proof}} \nc{\epf}{\end{proof}}
\nc{\be}{\begin{enumerate}}
\nc{\ee}{\end{enumerate}}
\nc{\bi}{\begin{itemize}}
\nc{\bimy}{\my{\begin{itemize}}
\nc{\eimy}{\end{itemize}}}
\nc{\itm}{\item}
\nc{\ei}{\end{itemize}}
\nc{\Subsection}[1]{\goodbreak\subsection*{#1}}%\ \par}

\nc{\sone}{\mathsf{S}_1}
\nc{\sfin}{\mathsf{S}_\mathrm{fin}}
\nc{\ufin}{\mathsf{U}_\mathrm{fin}}
\nc{\Split}{\mathsf{Split}}

\nc{\gone}{\mathsf{G}_1}    \nc{\gfin}{\mathsf{G}_\mathrm{fin}}

\DeclareMathOperator{\PK}{PK}

\DeclareMathOperator{\ds}{d}
\DeclareMathOperator{\ld}{ld}

\DeclareMathOperator{\wght}{w}
\DeclareMathOperator{\K}{K}
\DeclareMathOperator{\fin}{Fin}

\DeclareMathOperator{\C}{C}

\nc{\nbe}{\cN}
\nc{\ex}{\op{\exists}}
\nc{\Gbase}{$\mathfrak{G}$-base}
\nc{\Gkbase}{$\fin(\kappa)^\bbN$-base} %$\mathfrak{G}_\kappa$-base}
\nc{\Bd}[1]{\op{b}\Bigl(#1\Bigr)}
\nc{\bd}[1]{\op{b}(#1)}
\nc{\Bdd}[1]{\op{Bdd}(#1)}
\nc{\Bddd}{\op{Bdd}}
\nc{\bdd}{\mathrm{b}}
\nc{\Cof}[1]{\op{cof}\bigl(#1\bigr)}
\nc{\Cl}[2]{\op{cl}_{#1}(#2)}
\nc{\diam}{\op{diam}}
\nc{\Pp}{\mathrm{P}}
\nc{\IntN}{\bigcap_{n\in\bbN}}
\nc{\UnN}{\bigcup_{n\in\bbN}}
\nc{\cmpct}{\mathrm{c}}
\nc{\myb}[1]{{\color{blue}#1}}
\nc{\Nice}{\mathrm{N}}
\nc{\lbt}{local base trace}
\renewcommand{\poset}{partially ordered set}
\nc{\Poset}{Partially ordered set}
\nc{\grid}{grid}
\nc{\sgrid}{$\sigma$-grid}
\nc{\subgrid}{subgrid}
\nc{\ssubgrid}{$\sigma$-subgrid}

\newcommand{\seqn}[1]{{\la\, #1 \,\ra}}

\title[The cofinal structure of precompacta]{A classification of the cofinal structures of precompacta}
\keywords{Precompact sets, compact sets, compact-open topology, Tukey equivalence, local density}

\subjclass[2010]{
Primary:
03E15, %Descriptive set theory
03E04, %Ordered sets and their cofinalities; pcf theory
54A25; %Cardinality properties (cardinal functions and inequalities, discrete subsets)
Secondary:
03E17. %Cardinal characteristics of the continuum
}

\author[A. Eshed]{Aviv Eshed}
\address[Eshed]{Department of Mathematics, Weizmann Institute of Science, Rehovot, Israel}
\email{avivesh85@gmail.com}

\author[M. V. Ferrer]{M. Vicenta Ferrer}
\address[Ferrer]{Universitat Jaume I, IMAC and Departamento de Mate\-m\'{a}ticas,
Campus de Riu Sec, 12071 Castell\'{o}n, Spain.}
\email{mferrer@mat.uji.es}

\author[S. Hern\'andez]{Salvador Hern\'andez}
\address[Hern\'andez]{Universitat Jaume I, INIT and Departamento de Mate\-m\'{a}ticas,
Campus de Riu Sec, 12071 Castell\'{o}n, Spain.}
\email{hernande@mat.uji.es}

\author[P. Szewczak]{Piotr Szewczak}
\address[Szewczak]{
Institute of Mathematics, Faculty of Mathematics and Natural Science College of Sciences, Cardinal Stefan Wyszy\'nski University in Warsaw, Poland,
and
Department of Mathematics, Bar-Ilan University, Israel
}
\email{p.szewczak@wp.pl}
\urladdr{piotrszewczak.pl}

\author[B. Tsaban]{Boaz Tsaban}
\address[Tsaban]{Department of Mathematics, Bar-Ilan University, Israel}
\email{tsaban@math.biu.ac.il}
\urladdr{math.biu.ac.il/~tsaban}

\begin{document}

\begin{abstract}
We provide a complete classification of the possible cofinal structures of the families
of precompact (totally bounded) sets in general metric spaces,
and compact sets in general complete metric spaces.
Using this classification, we classify the cofinal structure of
local bases in the groups $\C(X,\bbR)$
of continuous real-valued functions on complete metric spaces $X$,
with respect to the compact-open topology.
\end{abstract}

\maketitle

\section{Introduction and related work}

A subset $C$ of a \poset{} $P$ is \emph{cofinal} if, for each element
$p\in P$, there is an element $c\in C$ with $p\leq c$. The cofinality of a \poset{} $P$, denoted by
$\cof(P)$, is the minimal cardinality of a cofinal subset of that set.
We identify the cofinality and, moreover, the cofinal structure of the family of precompact
sets in general metric spaces, and the family of compact sets in general complete metric spaces.
We apply our results to compute the character of the topological groups $\C(X,G)$
of continuous functions from a complete metric space $X$ to a group $G$ containing an arc,
equipped with the compact-open topology.
These extend results obtained earlier for metric \emph{groups}~\cite{PvK}.

The cofinal structure of the family of precompact sets in Polish spaces was identified by Christensen~\cite[Theorem~3.3]{christensen}.
Van Douwen~\cite{vD} computed the
cofinality of the \poset{} of \emph{compact} sets, and the
character of the diagonal sets~\cite[Theorem~8.13(c)]{vD}, for certain classes of separable metric spaces.
Using van Douwen's results, Nickolas and Tkachenko~\cite{NT1,NT2} computed
the character of the free topological groups $F(X)$
and free topological abelian groups $A(X)$ over the spaces $X$ considered
by van Douwen.
Independently of the present work,
and concurrently,
Gartside and Mamatelashvili~\cite{GM1,GM2,GMw1}
considered the cofinalities and the cofinal
structure of the partially ordered sets of
\emph{compact} sets for various spaces.
For non-scattered totally
imperfect separable metric spaces, and for
complete metric spaces of uncountable weight
less than $\aleph_\omega$, they identified the exact cofinal structures of families of compact sets.
The results of the paper cover all complete metric spaces.

For functions $f,g\in\NN$, define $f\leq g$ if $f(n)\leq g(n)$ for all natural numbers $n$.
For a topological group $G$,
let $\nbe(G)$ be the family of neighborhoods of the identity
element in $G$.
An \emph{$\NN$-base} (also known as \emph{\Gbase{}}) of a topological group $G$ is the image of a
a monotone cofinal map from $(\bbN^\bbN,\leq)$
to $(\nbe(G),\supseteq)$.
The notion of $\NN$-base has recently attracted considerable attention~\cite{CO,PvK,FK,GKL,GH,LPT,LRZ}.
%A generalization and a refinement of this notion,
%were, implicitly, studied in an earlier work~\cite{PvK}.
Gabriyelyan, K\k{a}kol, and Leiderman~\cite{GKL} apply Christensen's result to prove
that, for Polish spaces $X$, the topological groups $\C(X,\bbR)$
have $\NN$-bases.
Subsequently, Leiderman, Pestov and Tomita~\cite{LPT} proved that,
for collectionwise normal (e.g., paracompact) spaces $X$,
the character of the topological group $A(X)$ is determined by
the character of the diagonal set in the square $X\x X$.
A topological space is a \emph{$k$-space} if every set
with compact traces on all compact sets is closed.
Lin, Ravsky, and Zhang~\cite{LRZ} provide an inner characterization of
topological spaces such that the free topological group $F(X)$ is a
$k$-space with an $\NN$-base.
Banakh and Leiderman~\cite{BL} investigated various kinds of free objects of topological algebras with an $\NN$-base over separable, $\sigma$-compact, and metrizable spaces.

A $k$-space $X$ is a \emph{$k_\omega$-space} if it has a countable cofinal family
of compact sets.
The \emph{weight} of a topological space is the minimal cardinality of a basis of that space.
The \emph{compact weight} of a topological space is the supremum of the weights of its compact subsets.
Let $P$ and $Q$ be \poset{}s.
We write $P\preceq Q$
if there is a monotone cofinal map
$\psi\colon P\to Q$.
Two \poset{}s $P$ and $Q$ are \emph{cofinally equivalent}, denoted by $P\approx Q$, if $P\preceq Q$ and $Q\preceq P$.
For a set $X$, let $\fin(X)$ be the set of finite subsets of the set $X$. The family $\fin(X)^\bbN$ is the set of functions
$f \colon \bbN \to \fin(X)$.
For functions $f,g\in\fin(X)^\bbN$, define $f\leq g$ if $f(n)\sub g(n)$ for all natural numbers $n$.
Since the \poset{}s $\fin(\bbN)^\bbN$ and $\NN$ are cofinally equivalent,
a topological group $G$ has an
$\NN$-base if, and only if, $\fin(\bbN)^\bbN\preceq \nbe(G)$.
The following generalization and refinement of the notion of $\NN$-base were
considered, implicitly, earlier~\cite{PvK}.

\bdfn
Let $\kappa$ be an infinite cardinal number, and $G$ be a topological group.
A \emph{\Gkbase{}} of the group $G$ is the image of a monotone cofinal
function from the partially ordered set $\fin(\kappa)^\bbN$ into $\nbe(G)$.
\edfn

\bthm[{\cite[proof of Theorem 4.1]{PvK}}]
Let $X$ be a nondiscrete $k_\omega$-space of compact weight~$\kappa$, and $A(X)$ be the free abelian topological group over $X$.
Then $\cN(A(X))\approx \fin(\kappa)^\bbN$.
In particular, the group $A(X)$ has a $\fin(\kappa)^\bbN$-base.
\ethm

Let $X$ be a metric space.
Let $\PK(X)$ be the family of precompact sets in the space $X$.
The space $X$ is \emph{locally precompact} if every point in $X$ has a precompact neighborhood.

The \emph{density} of a topological space $X$, denoted by $\ds(X)$,
is the minimal cardinality of a dense subset.
The \emph{local density}~\cite{PvK} of a topological group $G$, denoted by $\ld(G)$,
is the minimal density of a neighborhood of the identity element (or, by homogeneity, any other group element).
Each item in the following theorem implies the subsequent one.

\bthm[{\cite[Corollary~3.7, Proposition~4.13, Theorem~4.16]{PvK}}]
Let $G$ be a metrizable group that is not locally precompact.
\be
\item Let $H$ be a clopen subgroup of the group $G$, of density $\ld(G)$.
Then $\PK(G)\approx\fin(G/H)\times \fin(\ld(G))^\bbN$.
\item If $\ld(G)=\ds(G)$,
then $\PK(G)\approx\fin(\ds(G))^\bbN$.
\item If the group $G$ is separable, then $\PK(G)\approx\NN$.
\ee
\ethm

A \poset{} $P$ is \emph{Tukey-finer} than a \poset{} $Q$~\cite{Tuk} if there is a map from the set $P$ into the set $Q$ carrying cofinal sets in $P$ into cofinal sets in $Q$.
\Poset{}s $P$ and $Q$ are \emph{Tukey-equivalent} if each of these sets is Tukey-finer than the other one.
The notion of cofinal equivalence used here is
stronger than Tukey-equivalence.
For our purposes, cofinal equivalence is also simpler.
Of course, Tukey-equivalence has its advantages.
Let $G$ be a topological group such that there is a cofinal family $\cB$ in $\cN(G)$ with $\cB\approx\fin(\kappa)^\bbN$.
This property may be insufficient for cofinal equivalence between $\cN(G)$ and $\fin(\kappa)^\bbN$ as defined here,
but it is sufficient for Tukey equivalence.
Readers more interested in Tukey equivalence may use this notion, instead
of cofinal equivalence, throughout the paper.

\section{Precompact sets}

\subsection{Preparations}

\bdfn
Let $X$ be a metric space, and $\epsilon$ be a positive real number.
A set $D\sub X$ is an $\epsilon$-\emph{\subgrid{}} if
the $\epsilon$-balls centered at the points of the set $D$ are
pairwise disjoint.
An $\epsilon$-\emph{\grid{}} is a maximal $\epsilon$-subgrid.
A set is a \emph{(sub)\grid{}} if it is an $\epsilon$-(sub)\grid{}
for some positive real number $\epsilon$.
\edfn

\blem\label{lem:fin}\mbox{}
\be
\item Every intersection of a precompact set and a sub\grid{} is finite.
\item Every nonprecompact metric space contains an infinite grid.
\ee
\elem

\bpf
(1) It follows from the definition of precompact sets.

(2) Let $X$ be a nonprecompact metric space.
Let $\epsilon$ be a positive real number such that the space $X$ cannot be covered by finitely many $\epsilon$-balls.
Let $D$ be an $\frac{\epsilon}{2}$-grid in $X$.
Then the family $\sset{B(x,\epsilon)}{x\in D}$ covers the space $X$.
Thus, the set $D$ is infinite.
\epf

\bdfn
A \emph{$\sigma$-(sub)grid} in a metric space $X$
is a family $D=\UnN D_n$ such that, for each natural number $n$, the set
$D_n$ is a $1/n$-(sub)grid containing the set $D_{n-1}$.
\edfn

\blem\label{lem:dense}
Every \ssubgrid{} in a metric space $X$ has cardinality at most $\ds(X)$.
In particular, every \sgrid{} in a metric space $X$ has cardinality $\ds(X)$.
\elem
\bpf
If $\ds(X)<\aleph_0$, then the space $X$ is finite.
Thus, we assume that $\ds(X)\geq\aleph_0$.
Let $D=\Un_{n\in\bbN}D_n$ be a $\sigma$-\subgrid{}.
For each natural number $n$, the set $D_n$ is a \subgrid{},
and thus $\card{D_n}\leq \ds(X)$.
It follows that
\[
\smallcard{D}= \medcard{\Un_{n\in\bbN}D_n}\leq \alephes\cdot \sup_{n\in\bbN}\card{D_n}\leq \aleph_0\cdot \ds(X)=\ds(X).
\]

If the set $D$ is a \sgrid{}, then it is dense, and thus $\card{D}=\ds(X)$.
\epf

\bprp\label{prp:dX<cfPK}
Let $X$ be a nonprecompact metric space. Then $\PK(X)\preceq \fin(\ds(X))$.
\eprp
\bpf
Let $D=\UnN D_n$ be a \sgrid{}.
By Lemma~\ref{lem:dense}, we have $\card{D}=\ds(X)$.
Thus, it suffices
to find a monotone cofinal map $\psi\colon\PK(X)\to \fin(D)$.

Since the metric space $X$ is not precompact, by Lemma~\ref{lem:fin}(2), there is a countably infinite subgrid
$\set{x_n}{n\in\bbN}$.
For each set $P\in\PK(X)$, let $n(P)$ be the maximal natural number such that $x_{n(P)}\in P$, or 1 if there is no such number.
Consider the map
\begin{align*}
\psi\colon\PK(X) &\longrightarrow \fin(D),\\
P &\longmapsto P\cap D_{n(P)}.
\end{align*}
By Lemma~\ref{lem:fin}(1), the map $\psi$ is well defined.
Since the sets $D_n$ increase with $n$, the
map $\psi$ is monotone.

The map $\psi$ is cofinal:
For a set $F\in\fin(D)$,
let $n$ be a natural number with $F\sub D_n$.
Then the set $P:=F\cup\{x_n\}$ is precompact, and $F\sub \psi(P)$.
\epf

\subsection{Metric spaces with a locally compact completion}

If the completion of a metric space is locally compact, then
the space is locally precompact. The converse implication fails.
Indeed, local compactness is not preserved by metric completions:
The hedgehog space $J(\alephes)$~\cite[Example~4.1.4]{eng}
is not locally compact, but
it is the completion of its locally compact subspace $J(\alephes)\sm\{0\}$.

\bthm\label{thm:SLPK}
Let $X$ be a nonprecompact metric space with a locally compact completion. Then $\PK(X)\approx \fin(\ds(X))$.
\ethm
\bpf
By Proposition~\ref{prp:dX<cfPK}, it remains to prove that $\fin(\ds(X)) \preceq \PK(X)$.

Let $D$ be a dense subset of the space $X$, of cardinality $\ds(X)$.
Let $\cB$ be the family of compact balls of rational radii
in the completion that are centered at $D$.
The closure of every precompact subset of $X$ in the completion is compact, and thus is covered by finitely many elements of the family $\cB$:
Fix a point $x\in X$.
There is a compact ball $B$ in the completion that is centered at the point $x$.
Since the set $D$ is dense in the completion, there is a ball of rational radii that is centered at some point of $D$, containing the point $x$, whose closure is contained in the ball $B$.

Since $\card{D}=\ds(X)$, we have $\card{\cB}=\ds(X)$, and
the union map,
\begin{align*}
\Un\colon \fin(\cB) &\longrightarrow \PK(X),\\
 \cF &\longmapsto \Un\cF\cap X,
\end{align*}

is monotone and cofinal.
\epf

\subsection{Metric spaces whose completion is not locally compact}

\bdfn
For an infinite cardinal number $\kappa$, let $\Bdd{\fin(\kappa)^\bbN}$ be the set
of the bounded functions in the set $\fin(\kappa)^\bbN$, that is,
$\Bdd{\fin(\kappa)^\bbN}:=\Un_{\alpha<\kappa}\fin(\alpha)^\bbN$.
\edfn

\blem\label{lem:Bddcf}
For cardinal numbers $\kappa$ of uncountable cofinality, we have  $\Bdd{\fin(\kappa)^\bbN}\allowbreak=\allowbreak\fin(\kappa)^\bbN$.\qed
\elem

\bdfn
Let $X$ be a topological space.
The \emph{local density} of a point $x\in X$, denoted by $\ld(x)$,
is the minimal density of a neighborhood of the point $x$.
The \emph{local density} of the space $X$, denoted by $\ld(X)$,
is the supremum of the local densities of the points of $X$.
The local density of the space $X$ is \emph{realized} if there is a point
$x\in X$ such that $\ld(x)=\ld(X)$.
\edfn

\bthm\label{thm:master}
Let $X$ be a metric space whose completion $\tilde X$ is not locally compact.
\be
\item If the local density of the completion is not realized, then
\[
\PK(X)\approx\fin(\ds(X))\x\Bdd{\fin(\ld(\tilde X))^\bbN}.
\]
\item If the local density of the completion is realized,
or its cofinality is uncountable,
then
\[
\PK(X)\approx\fin(\ds(X))\x\fin(\ld(\tilde X))^\bbN.
\]
\ee
\ethm

We prove Theorem~\ref{thm:master} using a series of Lemmata.

\blem\label{lem:cfPK<nd}
Let $X$ be a nonprecompact metric space with the completion $\tilde X$.
\be
\item $\fin(\ds(X))\x \fin(\ld(\tilde X))^\bbN\preceq\PK(X)$.

\item If the local density of the completion is not realized, then
\[
\fin(\ds(X))\x \Bdd{\fin(\ld(\tilde X))^\bbN}\preceq\PK(X).
\]
\ee
\elem
\bpf
(1) If the completion $\tilde{X}$ is locally compact, then the assertion follows from Theorem~\ref{thm:SLPK}.
Thus, assume that the completion is not locally compact.
The completion of the space $X$ with respect to an equivalent bounded metric is equal to $\tilde{X}$. Since density and local density do not depend on the metric, we assume that the metric of the space $X$ is bounded. For each point $x\in X$, let $r_x$ be the maximal radius such that $\ds(B(x,r_x))\leq\ld(\tilde{X})$.
Let $D$ be a \sgrid{} in $X$.
Fix a point $x\in X$.
The set $D\cap B(x,r_x)$ is a \ssubgrid{} in the set $B(x,r_x)$.
By Lemma~\ref{lem:dense}, we have $\card{D\cap B(x,r_x)}\leq\ld(\tilde{X})$.
Let $i_x\colon D\cap B(x,r_x)\to\ld(\tilde{X})$
be an injection.

For a subset $A$ of $X$, let $B(A,{1}/{n}):=\Un_{x\in A}B(x,{1}/{n})$.
Since $\card{D}=\ds(X)$, it suffices to prove that $\fin(D)\x \fin(\ld(\tilde X))^\bbN\preceq\PK(X)$.
Consider the map
\begin{align*}
\psi\colon \fin(D)\times \fin(\ld(\tilde{X}))^\mathbb{N} &\longrightarrow \PK(X),\\
(F,f) &\longmapsto \bigcap_{n\in \mathbb{N}} B\Bigl(\Un_{x\in F}i_x\inv[f(n)],{1}/{n}\Bigr).
\end{align*}
Since the sets $\Un_{x\in F}i_x\inv[f(n)]$ are finite for all natural numbers $n$, the set $\psi(F,f)$ is precompact. Thus, the map $\psi$ is well defined.
By the definition, the map $\psi$ is monotone.

Every precompact set in the space $X$ is covered by finitely many balls $B(x,r_x)$ centered at the points of $D$:
For a ball $B(x,r)$ in the space $X$, let $\tilde{B}(x,r)$ be the corresponding $r$-ball in the completion.
Since the closure of each precompact set in the completion is compact, it is enough to show that the family $\sset{\tilde{B}(x,r_x)}{x\in D}$ is a cover of $\tilde{X}$:
Fix a point $\tilde{x}\in\tilde{X}$.
Let $U$ be an open neighborhood of $\tilde{x}$, of density $\ld(\tilde{x})$.
There are a point $x\in D\cap U$, and a positive real number $r$ such that $\tilde{x}\in\tilde{B}(x,r)\sub U$.
Then $\ds(B(x,r))=\ld(\tilde{x})\leq\ld(\tilde X)$. Thus, $r\leq r_x$.

The map $\psi$ is cofinal:
Let $P$ be a precompact set in $X$. By the above, there is a finite set $F\sub D$ such that $P\sub \Un_{x\in F}B(x,r_x)$.
For each natural number $n$, there is a finite set $L_n\sub \Un_{x\in F}B(x,r_x)\cap D$ such that $P\sub B(L_n,{1}/{n})$.
Define a function $f\in \fin(\ld(\tilde{X}))^\mathbb{N}$ by
 $f(n):=\Un_{x\in F}i_x(L_n\cap B(x,r_x))$ for all natural numbers $n$. Fix a natural number $n$.
Since $L_n\sub \Un_{x\in F}i_x\inv[f(n)]$, we have
\[
P\sub B\bigl(L_n,{1}/{n}\bigr)\sub B\bigl( \Un_{x\in F}i_x\inv[f(n)],{1}/{n}\bigr).
\]

Thus,
\[
P\sub \bigcap_{n\in \mathbb{N}} B\bigl(\Un_{x\in F}i_x\inv[f(n)],{1}/{n}\bigr)=\psi(F,f).
\]

(2) We make the following modifications in the proof of the first item:
Fix a point $x\in X$.
Let
$r_x$ be the supremum of the radii $r$ such that $\ds(B(x,2r))<\ld(\tilde{X})$.
Then $\ds(B(x,r_x))<\ld(\tilde{X})$.
Let $\psi$ be the map defined in~(1).
Let $P$ be a precompact set in the space $X$,
$F\sub D$ be a finite set such that $P\sub\Un_{x\in F}B(x,r_x)$, and $\kappa:=\max\sset{\im[i_x]}{x\in F}$.
Then $\kappa<\ld(\tilde{X})$, and the function $f\in\fin(\ld(\tilde{X}))^\bbN$, defined in the proof of~(1), is bounded by $\kappa$.
\epf

\blem\label{lem:sup<cfK}
Let $X$ be a metric space, and $x$ be a point in the completion with no compact neighborhood.
Then $\PK(X)\preceq\fin(\ld(x))^{\bbN}$.
\elem
\bpf
Assume first that the local density of the completion is countable.
If $\ld(x)$ is finite, then $x$ is an isolated point in the completion, and then the point $x$ has a compact neighborhood, a contradiction.
Thus, $\ld(x)=\aleph_0$.
Let $\sset{U_n}{n\in\bbN}$ be a decreasing local base in the completion at the point $x$  such that $\ds(U_n)=\aleph_0$ for all natural numbers $n$.
For each natural number $n$, the set $U_n\cap X$ is nonprecompact, and by Lemma~\ref{lem:fin}(2),
there is a countably infinite grid $D_n$ in the set $U_n\cap X$.
Thus, $\prod_n\fin(D_n)\approx \fin(\aleph_0)^\bbN$.

Define a monotone map
\begin{align*}
\psi\colon \PK(X) &\longrightarrow \prod_{n\in\bbN}\fin(D_n),\\
P &\longmapsto \seq{P\cap D_n}{n\in\bbN}.
\end{align*}
By Lemma~\ref{lem:fin}(1), the map $\psi$ is well defined.
Fix a function $f\in \prod_{n\in\bbN}\fin(D_n)$.
Let $P:=\Un_{n\in\bbN} f(n)$.
For each positive real number $\epsilon$, the set $P\sm B(x,\epsilon)$ is finite.
Thus, the set $P$ is precompact, and $f\leq\psi(P)$.
The map $\psi$ is cofinal.
We have $\PK(X)\preceq \prod_n\fin(D_n)\approx\fin(\aleph_0)^\bbN$.

Next, assume that the local density of the completion is uncountable.
Pick a point $y$ in the completion with $\ld(y)\geq\max\{\ld(x),\aleph_1\}$.
Since $\fin(\ld(y))^\bbN\preceq \fin(\ld(x))^\bbN$, it suffices to prove the assertion for the point $y$.
Let $\kappa:=\ld(y)$.
Let $\sset{U_n}{n\in\bbN}$ be a decreasing local base in the completion at the point $y$  such that $\ds(U_n)=\kappa$ for all natural numbers $n$.
Let $D=\Un_{n\in\bbN}D_n$ be a \sgrid{} in the space $X$.
Fix a natural number $n$.
The set $D\cap U_n$ is a dense $\sigma$-subgrid of the set $U_n$.
By Lemma~\ref{lem:dense}, we have $\card{D\cap U_n}=\kappa$.
The sets $D_m$ increase with the index $m$, and $\sup_m\card{D_m\cap U_n}=\kappa$.

There is a sequence $\la m_n \ra_{n=1}^\infty$ of natural numbers such that the cardinal numbers $\kappa_n:=\card{D_{m_n}\cap U_n}$ are nondecreasing with the index $n$, and $\sup_n\kappa_n=\kappa$:
If there is a sequence $\la m_n \ra_{n=1}^\infty$ of natural numbers such that $\card{D_{m_n}\cap U_n}=\kappa$ for all natural numbers $n$, take $\kappa_n:=\card{D_{m_n}\cap U_n}$.
Thus, assume that there is no such sequence.
Since the sets $U_n$ decrease with the index $n$, we have $\card{D_m\cap U_n}<\kappa$ for all but finitely many natural numbers $n$, and all natural numbers $m$.
Removing the first few sets $U_n$, we assume that the inequality $\card{D_m\cap U_n}<\kappa$ holds for all natural numbers $n$, and $m$.
For each natural number $n$, there is a natural number $m_n$ such that the cardinal number $\kappa_n:=\card{D_{m_n}\cap U_n}$ is greater than the cardinal numbers $\card{D_n\cap U_1}$, and $\kappa_{n-1}$.

Since the cardinal numbers $\kappa_n$ are nondecreasing with the index $n$, we have
\[
\prod_{n\in\bbN}\fin(D_{m_n}\cap U_n)\approx\prod_{n\in\bbN}\kappa_n\approx\fin(\kappa)^\bbN\text{~\cite[Lemma~4.14]{PvK}}.
\]

Define a monotone map
\begin{align*}
\psi\colon \PK(X) &\longrightarrow \prod_{n\in\bbN}\fin(D_{m_n}\cap U_n),\\
P &\longmapsto \seq{P\cap D_{m_n}\cap U_n}{n\in\bbN}.
\end{align*}
As in the previous case, the sets $D_{m_n}\cap U_n$ are subgrids for all natural numbers $n$, and thus the map $\psi$ is well defined and cofinal.
Then $\PK(X)\preceq \prod_n\fin(D_{m_n}\cap U_n)\approx\fin(\kappa)^\bbN$.
\epf

\blem\label{lem:Bdd}
Let $X$ be a metric space whose completion is not locally compact.
Let $\kappa$ be the local density of the completion of the space $X$, $\lambda$ be its cofinality, and
$\seqn{\kappa_\alpha}_{\alpha<\lambda}$ be an increasing sequence of cardinal numbers with supremum $\kappa$.
If $\PK(X)\preceq \fin(\kappa_\alpha)^\bbN$ for each ordinal number $\alpha<\lambda$, then $\PK(X)\preceq\Bdd{\fin(\kappa)^\bbN}$.
\elem
\bpf
The proof generalizes the proof of Proposition~\ref{prp:dX<cfPK}.
There is a \subgrid{} of cardinality $\lambda$:
Let $D=\UnN D_n$ be a \sgrid{}. Since $X$ is nonprecompact, by Lemma~\ref{lem:fin}(2), there is a countably infinite \subgrid{}.
Thus, if $\lambda$ is countable, we are done.
Assume that $\lambda$ is uncountable.
If $\lambda<\ds(X)$, then by Lemma~\ref{lem:dense},  there is a natural number $k$ such that $\card{D_k}>\lambda$.
If $\lambda=\ds(X)$, then $\lambda\geq \ld(\tilde X)=\kappa$, which implies that the cardinal number $\kappa$ is regular.
Then there is a natural number $k$ such that $\card{D_k}=\ds(X)$.

Let $A=\set{a_\alpha}{\alpha<\lambda}$ be a \subgrid{}. For each ordinal number $\alpha<\lambda$,
let $\psi_\alpha\colon\PK(X)\to\fin(\kappa_\alpha)^\bbN$ be a monotone cofinal map.
By Lemma~\ref{lem:fin}(1), for each precompact set $P$ in $X$, there is a finite set $F(P)\sub\lambda$ such that $P\cap A=\set{a_\alpha}{\alpha\in F(P)}$.
Define a map
\begin{align*}
\psi\colon\PK(X) &\longrightarrow\Bdd{\fin(\kappa)^\bbN},\\
P &\longmapsto \psi(P)(n):=\Un_{\alpha\in F(P)}\psi_\alpha(P)(n).
\end{align*}

The map $\psi$ is monotone: Fix precompact sets $P$ and $P'$ in $X$ such that $P\sub P'$, and a natural number~$n$. Then $F(P)\sub F(P')$, and
\[
\psi(P)(n)=\Un_{\alpha\in F(P)} \psi_\alpha(P)(n)\sub \Un_{\alpha\in F(P')} \psi_\alpha(P)(n)\sub \Un_{\alpha\in F(P')} \psi_\alpha(P')(n)=\psi(P')(n).
\]
Thus, $\psi(P)\leq\psi(P')$.

The map $\psi$ is cofinal: Fix a function $f\in\Bdd{\fin(\kappa)^\bbN}$.
There is an ordinal number $\beta<\lambda$ such that $f\in\fin(\kappa_\beta)^\bbN$.
Then there is a precompact set $P$ in $X$ such that $f\leq\psi_\beta (P)$.
The set $P':=P\cup\{a_\beta\}$ is precompact in $X$. For each natural number $n$, we have
\[
f(n)\sub \psi_\beta(P)(n)\sub\psi_\beta(P')(n)\sub\Un_{\alpha\in F(P')}\psi_\alpha(P')(n)=\psi(P')(n).
\]
Thus, $f\leq\psi(P')$.
\epf

\bpf[Proof of Theorem~\ref{thm:master}]
(1)
Let $\lambda$ be the cofinality of $\ld(\tilde{X})$.
Since the cardinal number $\ld(\tilde{X})$ is uncountable, and each compact metric space is separable, there is a set $\sset{x_\alpha}{\alpha<\lambda}$ of points of the completion with no compact neighborhoods such that the sequence $\la \ld(x_\alpha)\ra_{\alpha<\lambda}$ is increasing, with supremum $\ld(\tilde{X})$.
By Lemmata~\ref{lem:sup<cfK} and~\ref{lem:Bdd}, we have $\PK(X)\preceq\Bdd{\fin(\ld(\tilde{X}))^\bbN}$.
By Proposition~\ref{prp:dX<cfPK}, we have
$\PK(X)\preceq \fin(\ds(X))\x\Bdd{\fin(\ld(\tilde X))^\bbN}$.
Apply Lemma~\ref{lem:cfPK<nd}(2).

(2) Assume that the local density of the completion is realized.
There is a point $x$ in the completion with $\ld(x)=\ld(\tilde{X})$.
By Proposition~\ref{prp:dX<cfPK} and Lemma~\ref{lem:sup<cfK}, we have
$\PK(X)\preceq\fin(\ds(X))\x\fin(\ld(\tilde X))^\bbN$.
Apply Lemma~\ref{lem:cfPK<nd}(1).
If the local density of the completion is not realized,
and its cofinality is uncountable, apply~(1) and Lemma~\ref{lem:Bddcf}.
\epf

\subsection{Inner characterizations}

The results of this section can be stated in an inner manner, that is,
directly in terms of properties of the metric space $X$, without consideration of its
completion. This is done in the following proposition, whose straightforward
proofs are omitted.

\bdfn
Let $X$ be a metric space.
A \emph{\lbt{}} in $X$ is a family
of open sets whose members are the traces, in $X$, of the elements
of a local base in the completion of the space $X$.
\edfn

\bprp
Let $X$ be a metric space.
\be
\item
A family of open sets in the space $X$ is a \lbt{} if and only if
every finite intersection of members of this family contains
the closure of some member of the family, and the diameters of the sets
in the family are not bounded away from 0.
\item The local density of the completion of the space $X$ is
equal to the supremum of the cardinal numbers
$\min\sset{\ds(U)}{U\in\cU}$, for the \lbt{}s $\cU$ in the space $X$.
\item The completion of a metric space is locally compact if and only
if every \lbt{} in that space contains a precompact set.
\item The local density of the completion of a metric space $X$ is realized
if and only if the local density of the completion is equal to
the cardinal number $\min\sset{\ds(U)}{U\in\cU}$, for some
\lbt{} $\cU$ in the space $X$.\qed
\ee
\eprp

\section{Compact sets}

For a topological space $X$, let $\K(X)$ be the family of compact subsets of the space~$X$.

\bprp\label{prp:dX<cfK}
Let $X$ be a noncompact metric space. Then $\K(X)\preceq \fin(\ds(X))$.
\eprp
\bpf
There is a countably infinite closed discrete set in the space $X$.
Proceed as in the proof of Proposition~\ref{prp:dX<cfPK}.
\epf

\bthm\label{thm:masterComplete}
Let $X$ be a complete metric space.
\be
\item If the space $X$ is noncompact and locally compact, then $\K(X)\approx\fin(\ds(X))$.
\item
If the space $X$ is not locally compact, then:
\be
\item If the local density of the space $X$ is not realized,
then $\K(X)\approx\fin(\ds(X))\x\Bdd{\fin(\ld(X))^\bbN}$.
\item If the local density of the space $X$ is realized,
or its cofinality is uncountable, then
$\K(X)\approx\fin(\ds(X))\x\fin(\ld(X))^\bbN$.
\ee
\ee
\ethm
\bpf
(1) By Proposition~\ref{prp:dX<cfK}, it suffices to prove that $\fin(\ds(X)) \preceq \K(X)$.
Proceed as in the proof of Theorem~\ref{thm:SLPK},
with the exception that the family $\cB$ is an open locally finite cover of the space $X$ by sets whose closures are compact~\cite[Theorem~4.4.1]{eng}.

(2) For a complete metric space $X$, we have $\K(X)\approx \PK(X)$:
From the set $\PK(X)$ to the set $\K(X)$, we take the map $P\mapsto \cl{P}$.
In the other direction, we take the identity map.
Apply Theorem~\ref{thm:master}.
\epf

\bcor\label{cor:K=}
Let $X$ be a complete metric space that is not locally compact,
with $\ds(X)=\ld(X)$.
If the local density of the space $X$ is realized, or its cofinality
is uncountable, then
$\K(X)\approx\fin(\ds(X))^\bbN$.\qed
\ecor

Let $X$, $Y$ be topological spaces, $K$ be a compact set in the space $X$, and $U$ be an open set in the space $Y$.
Let $[K,U]$ be the set of functions $f\in \C(X,Y)$ such that $f[K]\sub U$, a basic open set in $\C(X,Y)$

\blem\label{lem:chiCXG}
Let $X$ be a Tychonoff space, and $G$ be a topological group.
If the sets $\cK$ and $\cU$ are cofinal in the sets $\K(X)$ and $\nbe(G)$, respectively, then $\cK\x\cU\preceq \nbe(\C(X,G))$.
\elem

\bpf
The map $(K,U)\mapsto [K,U]$, from the set $\cK\x\cU$ into the set $\cN(\C(X,G))$, is monotone and cofinal.
\epf

\bprp\label{prp:chiCXG}
Let $X$ be a nonempty Tychonoff space, and $G$ be a topological group containing an arc.
Let $\cK$ and $\cU$ be directed subsets of the sets $\K(X)$ and $\nbe(G)$, respectively.
If the set $\cB:=\sset{[K,U]}{K\in\cK,U\in\cU}$ is cofinal in the set $\nbe(\C(X,G))$, then $\cB\approx\cK\x\cU$, and the sets $\cK$ and $\cU$ are cofinal in the sets $\K(X)$ and $\nbe(G)$, respectively.
\eprp

In order to prove Proposition~\ref{prp:chiCXG}, we need the following results.

\blem\label{lem:KU}
Let $X$ be a Tychonoff space, and $Y$ be a topological space.
Let $K$, and $C$ be compact sets in the space $X$ with $C\neq\emptyset$, and $U$,and $V$ be open sets in the space $Y$ such that $[K,U]\sub[C,V]$.
Then:
\be
\item $U\sub V$.
\item If the space $Y$ contains an arc $A$
such that $A\cap U\neq\emptyset$ and $A\sm V\neq\emptyset$, then $C\sub K$.
\ee
\elem

\bpf
(1)
For each element $y\in U$, the constant function $c_y$ with value $y$ belongs to the set $[K,U]$.
Thus, $c_y\in [C,V]$.
For a point $x\in C$, we have $y=c_y(x)\in V$.
Thus, $U\sub V$.

(2)
Let $a\in A\cap U$, and $b\in A\sm V$.
Suppose that there is an element $x\in C\sm K$.
Since the set $A$ is homeomorphic to the unit interval, there is a function $f\in\C(X,Y)$ such that $f[K]=\{a\}$ and $f(x)=b$.
Since $a\in U$ and $b\nin V$, we have $f\in [K,U]\sm [C,V]$, a contradiction.
\epf

\blem\label{lem:poset}
Let $(P,\leq)$ be a directed set, and $q\in P$.
Then $P\approx\sset{p\in P}{q\leq p}$, the cone of all elements above the element $q$.
\elem

\bpf
The map $\psi\colon P\to \sset{p\in P}{q\leq p}$ defined by
\[
\psi(p):=
\begin{cases}
p,&\text{ if }q\leq p,\\
q,&\text{ otherwise,}
\end{cases}
\]
for all elements $p\in P$, is monotone and cofinal.
In the other direction, take the identity map.
\epf

\begin{proof}[Proof of Proposition~\ref{prp:chiCXG}]

From the set $\cK\x\cU$ to the set $\cB$, take the map $(K,U)\mapsto [K,U]$.
Thus, $\cK\x\cU\preceq\cB$.

Take an arc $A$ in the group $G$, containing the identity element $1_G$ of the group $G$.
Since the set $\cB$ is cofinal in $\nbe(\C(X,G))$, there are a nonempty set $L\in\cK$, and a set $W\in\cU$ with $A\not\sub U$.
By Lemma~\ref{lem:poset}, we have
\[
\cK\approx \sset{K\in\cK}{L\sub K},\text{ and }\cU\approx\sset{U\in\cU}{U\sub W}.
\]
By Lemmata~\ref{lem:KU} and~\ref{lem:poset}, we have
\[
\cB\approx\sset{[K,U]}{K\in\cK, L\sub K,U\in\cU, U\sub W}.
\]
Thus assume that the sets $\cK$, $\cU$, and $\cB$ are equal to the cones just defined.

Define a map $\psi\colon\cB\to \K(X)\x\nbe(G)$ by $\psi([K,U]):=(K,U)$ for all sets $[K,U]\in\cB$.
By Lemma~\ref{lem:KU},
the map $\psi$ is well defined and monotone.
The map $\psi$ is cofinal:
Let $(C,V)\in \K(X)\x\nbe(G)$.
Take an element $a\in A\cap V$ such that $a\neq 1_G$.
Since the set $\cB$ is cofinal in $\nbe(\C(X,G))$, there is a set $[K,U]\in \cB$ such that $[K,U]\sub [C, V\sm\{a\}]$.
By Lemma~\ref{lem:KU}, we have $C\sub K$ and $U\sub V$.
The set $\cK\x \cU$ is the image of the map $\psi$,
which is cofinal in the set $\K(X)\x\nbe(G)$.
Thus, $\cB\preceq \cK\x \cU$, and the sets $\cK$ and $\cU$ are cofinal in the sets $\K(X)$ and $\nbe(G)$, respectively.
\epf

\bthm\label{thm:exactGbase}
Let $X$ be a complete metric space that is not locally compact,
with $\ds(X)=\ld(X)$.
Let $G$ be a metrizable topological group.
Assume that the local
density of the space $X$ is realized, or its cofinality
is uncountable.

\be
\item The topological group $\C(X,G)$ has a $\fin(\ds(X))^\bbN$-base.
\item If the group $G$ contains an arc, then the group $\C(X,G)$ has a local base at the identity element that is cofinally equivalent to $\fin(\ds(X))^\bbN$
\ee
\ethm

\bpf
(1)
Let $\kappa:=\ds(X)$.
By Corollary~\ref{cor:K=}, we have $\K(X)\approx\fin(\kappa)^\bbN$.
Since the group $G$ is metrizable, we have $\nbe(G)\approx\fin(\w)$.
By Lemma~\ref{lem:chiCXG}, we have
\[
\fin(\kappa)^\bbN\approx\fin(\kappa)^\bbN\x \fin(\w)\preceq\nbe(\C(X,G)).
\]

(2)
Let $\kappa:=\ds(X)$, and  $\cB:=\sset{[K,U]}{K\in\K(X),U\in\nbe(G)}$.
By Proposition~\ref{prp:chiCXG}, we have $\fin(\kappa)^\bbN\approx \fin(\kappa)^\bbN\x\fin(\w)\approx \cB$, and $\cB\preceq\nbe(\C(X,G))$.
\epf

\brem
Theorem~\ref{thm:exactGbase} remains true if we consider a topological group $G$ with a cofinal set $\cU$ in $\nbe(G)$ such that $\cU\approx\fin(\lambda)$ or $\cU\approx\fin(\lambda)^\bbN$ for some cardinal number $\lambda\leq \ds(X)$.
\erem
\section{Cofinality}

Theorem~\ref{thm:SLPK} implies the following result.

\bcor
Let $X$ be a nonprecompact metric space with a locally compact completion.
Then $\cof(\PK(X))=\ds(X)$.\qed
\ecor

When the completion is not locally compact, we need to estimate the
cofinality of the partially ordered set $\fin(\kappa)^\bbN$. This cofinality
can be expressed by well-studied set theoretic functions~\cite[Section~8]{PvK}.
The introductory material in this section is adapted from an earlier paper~\cite{Psi}.
For an infinite cardinal number $\kappa$, let $[\kappa]^\alephes$
be the family of all countably infinite subsets of $\kappa$.
Let $\fc$ be the cardinality of the continuum.
Let $\fd:=\cof(\NN,\le)$. We have $\aleph_1\le\fd\le\fc$, and the cofinality of
the cardinal number $\fd$ is uncountable~\cite{BlassHBK}.

\bprp[{\cite[Proposition~4.15]{PvK}}]\label{prp:morph}
Let $\kappa$ be an infinite cardinal number. Then
$\cof(\fin(\kappa)^\bbN)=\fd\cdot\cof([\kappa]^\alephes)$.
\eprp

The estimation of the cardinal number $\cof([\kappa]^\alephes)$ in terms of the cardinal number $\kappa$
is a central goal in Shelah's \emph{PCF theory}, the theory of possible cofinalities.
The \emph{PCF function} $\kappa\mapsto\cof([\kappa]^\alephes)$ is tame.
For example, if there are no large cardinals in the Dodd--Jensen core model,
then $\cof([\kappa]^\alephes)$ is simply $\kappa$ if $\kappa$ has uncountable cofinality,
and $\kappa^+$ (the successor of $\kappa$) otherwise.
Moreover, without any special hypotheses, the cardinal number
$\cof([\kappa]^\alephes)$ can be estimated, and in many cases computed exactly.
Some examples follow.

For uncountable cardinal numbers $\kappa$ of countable cofinality,
a variation of K\"onig's Lemma implies that $\cof([\kappa]^\alephes)>\kappa$.
\emph{Shelah's Strong Hypothesis (SSH)} is the assertion that
$\cof([\kappa]^\alephes)=\kappa^+$ for all uncountable cardinal numbers $\kappa$ of countable cofinality.
The Generalized Continuum Hypothesis implies SSH, but the latter axiom is much
weaker, being a consequence of the absence of large cardinal numbers.

\bthm[Folklore]\label{nice}
The following cardinal numbers are fixed points of the PCF function $\kappa\mapsto\cof\allowbreak([\kappa]^\alephes)$:
\be
\item The cardinal numbers $\kappa$ with $\kappa^\alephes=\kappa$.
\item The cardinal numbers $\aleph_n$, for all natural numbers $n\ge 1$.
\item The cardinal numbers $\aleph_\kappa$, for a singular cardinal number $\kappa$ of uncountable cofinality
that is smaller than the first fixed point of the $\aleph$ function.
\item Assuming SSH, all cardinal numbers of uncountable cofinality.
\ee
Moreover, successors of fixed points of the PCF function are also fixed points of that function.
\ethm

For example, the cardinal numbers $\aleph_{\aleph_{\w_n}}$ and its successors are all fixed
points of the PCF function for all natural numbers $n$.
The following theorem summarizes some known anomalies of the PCF function,
and the function $\kappa\mapsto \cof(\fin(\kappa)^\bbN)$, in light of Proposition~\ref{prp:morph}.

\bthm[{\cite[\S 3]{Psi}}]
\mbox{}
\be
\item Let $\aleph_\alpha:=\fd$. If $\aleph_{\alpha+\omega}<\fc$,
then there is a cardinal number $\kappa<\fc$ such that $\cof(\fin(\kappa)^\bbN)>\fd\cdot\kappa$.
\item Assume SSH. Then:
\be
\itm For each infinite cardinal number $\kappa\le\fd$, we have $\cof(\fin(\kappa)^\bbN)=\fd$.
\itm We have $\cof(\fin(\kappa)^\bbN)=\fd\cdot\kappa$ for all infinite cardinal numbers $\kappa\le\fc$ if, and only if,
there is an integer $n\ge 0$ such that $\fc=\fd^{+n}$, the $n$-th successor of $\fd$.
\ee
\item It is consistent (relative to the consistency of ZFC with appropriate large cardinal numbers hypotheses) that
\[
\aleph_\w<\fd=\aleph_{\w+1}<\cof(\fin(\aleph_\w)^\bbN)=
\cof([\aleph_\w]^\alephes)=\aleph_{\w+\gamma+1}=\fc,
\]
for each prescribed ordinal number $\gamma$ with $1\le\gamma<\aleph_1$.
\ee
\ethm

\bthm\label{thm:cof}
Let $X$ be a metric space whose completion $\tilde X$ is not locally compact.
Then $\Cof{\PK(X)}=\ds(X)\cdot\fd\cdot\sup_{\tilde x\in\tilde X} \cof([\ld(\tilde x)]^\alephes)$.
\ethm
\bpf
If the local density of the completion is realized, apply Theorem~\ref{thm:master}(2).
Thus, assume that the local density of the completion is not realized.

\blem%\label{lem:Bddsup}
Let $\kappa$ be a cardinal number. Then
\[
\cof(\Bdd{\fin(\kappa)^\bbN})=\kappa\cdot\sup_{\alpha<\kappa}\Cof{\fin(\alpha)^\bbN}.
\]
\elem
\bpf
For each ordinal number $\alpha<\kappa$, let $Y_\alpha$ be a cofinal set in $\fin(\alpha)^\bbN$ of cardinality
$\Cof{\fin(\alpha)^\bbN}$.
The set $\Un_{\alpha<\kappa}Y_\alpha$ is cofinal in $\Bdd{\fin(\kappa)^\bbN}$, and thus
\[
\Cof{\Bdd{\fin(\kappa)^\bbN}}\leq\medcard{\Un_{\alpha<\kappa}Y_\alpha}\leq\kappa\cdot\sup_{\alpha<\kappa}\card{Y_\alpha}.
\]
Let $Y$ be a cofinal set in $\Bdd{\fin(\kappa)^\bbN}$ of cardinality $\Cof{\Bdd{\fin(\kappa)^\bbN}}$.
Fix an ordinal number $\alpha<\kappa$. For each function $f\in Y$,
define a function $f_\alpha\in \fin(\alpha)^\bbN$ by $f_\alpha(n):=f(n)\cap\alpha$ for all natural numbers $n$. The set $\set{f_\alpha}{f\in Y}$ is cofinal in $\Bdd{\fin(\kappa)^\bbN}$.
Since $\kappa\leq\card{Y}$, we have
\[
\kappa\cdot\sup_{\alpha<\kappa}\Cof{\fin(\alpha)^\bbN}\leq\card{Y}.\qedhere
\]
\epf

By Theorem~\ref{thm:master}(1), we have
\[\Cof{\PK(X)}=\ds(X)\cdot\ld(\tilde{X})\cdot\fd\cdot\sup_{\tilde x\in\tilde X} \cof([\ld(\tilde x)]^\alephes).
\]
Since $\ld(\tilde{X})\leq\ds(\tilde{X})=\ds(X)$, we are done.
\epf

By Theorem~\ref{thm:master} and Lemma~\ref{lem:Bdd}, we have the following corollary.

\bcor
Let $X$ be a metric space whose completion $\tilde X$ is not locally compact.
If the local density of the completion $\tilde X$ is realized,
or its cofinality is uncountable, then
\[
\cof(\PK(X))=\ds(X)\cdot\fd\cdot\cof\bigl([\ld(\tilde X)]^\alephes\bigr).\qed
\]
\ecor

\bthm\label{thm:hi5}
Let $X$ be a complete metric space.
\be
\item If the space $X$ is noncompact and locally compact, then $\cof(\K(X))=\ds(X)$.
\item
If the space $X$  is not locally compact, then
\[
\cof(\K(X))=\ds(X)\cdot\fd\cdot\sup_{x\in X}\Cof{[\ld(x)]^\alephes}.
\]
\ee
\ethm
\bpf
(1) Apply Theorem~\ref{thm:masterComplete}.

(2) Apply Theorem~\ref{thm:cof}.
\epf

\bcor\label{cor:masterComplete}
Let $X$ be a complete metric space that is not locally compact.
Assume that the local density of the space $X$ is realized,
or has uncountable cofinality. Then:
\be
\item
$\cof(\K(X))=\ds(X)\cdot\fd\cdot\Cof{[\ld(X)]^\alephes}$.
\item
If $\ds(X)=\ld(X)$, then
$\cof(\K(X))=\fd\cdot\cof([\ds(X)]^\alephes)$.
\ee
\ecor
\bpf
(1) Apply Theorem~\ref{thm:hi5}(2).

(2) Apply Corollary~\ref{cor:K=}.
\epf

\bexm
The equality from Corollary~\ref{cor:masterComplete}(1) is not provable for general metric spaces.
Assume \CH{}.
For a cardinal number $\kappa$, let $J(\kappa)$ be the corresponding \emph{hedgehog space}~\cite[Example~4.1.4]{eng}.
The space $X:=\bigoplus_{n\in\bbN} J(\aleph_n)$ is complete, and it is not locally compact. The local density $\aleph_\omega$ of $X$ is not realized.
Since $\ds(X)=\ld(X)=\aleph_\omega$, by Theorem~\ref{thm:masterComplete}(2a) and an earlier result~\cite[Proposition~4.15]{PvK}, we have
\[
\cof(\K(X))=\aleph_\omega\cdot\fd\cdot \sup_{n\in\bbN}\Cof{[\aleph_n]^\alephes}=\aleph_\omega\cdot\aleph_1\cdot\sup_{n\in\bbN}\aleph_n=\aleph_\omega\cdot\aleph_\omega=\aleph_\omega.
\]
On the other hand, we have
\[
\ds(X)\cdot\fd\cdot \Cof{[\ld(X)]^\alephes}=\aleph_\omega\cdot\aleph_1\cdot \Cof{[\aleph_\omega]^\alephes}=\aleph_\omega\cdot\aleph_\omega^+=\aleph_\omega^+.\qedhere
\]
\eexm

\section{Weight}

\bthm\label{thm:LieGen}
Let $X$ be a complete nondiscrete noncompact metric space, and $G$
be a topological group containing an arc.
\be
\item If the space $X$ is locally compact, then
\[
\ds(X)\leq\wght(\C(X,G))\leq\wght(G)\cdot\ds(X).
\]
\item If the space $X$ is not locally compact, then
\[
\ds(X)\cdot\fd\cdot\sup_{x\in X}\Cof{[\ld(x)]^\alephes}\leq \wght(\C(X,G))\leq
\wght(G)\cdot\ds(X)\cdot\fd\cdot\sup_{x\in X}\Cof{[\ld(x)]^\alephes}.
\]
\ee
\ethm

We prove Theorem~\ref{thm:LieGen} using the following notions and observations.
The \emph{character} $\chi(G)$ of a topological group $G$ is the minimal cardinality of a local
base in the group $G$.
The boundedness number $\bd{G}$ of a topological group $G$ is the minimal cardinal number such
that, for each open set $U$ in $G$, there is a set $S\subseteq G$ of cardinality $\bd{G}$
with $S\cdot U=G$.

\blem\label{lem:bdProd}
Let $\prod_{t\in T}G_t$ be a Tychonoff product of topological groups.
Then \\
$
\bd{\prod_{t\in T}G_t}=\sup_{t\in T}\bd{G_t}.
$
\elem

\bpf
Let $U$ be a basic open neighborhood of the identity element in the group $\prod_{t\in T}G_t$.
For each index $t\in T$, let $1_t$ be the identity element in the group $G_t$.
There are a finite set $F\sub T$,  and for each index $t\in T$, an open neighborhood $U_t$ of the element $1_t$ such that $U=\prod_{t\in F}U_t\x\prod_{t\in T\sm F}G_t$.
For each index $t\in T$, let $S_t$ be a subset of the group $G_t$ of cardinality $\bd{G_t}$  such that $S_t\cdot U_t=G_t$.
We have $(\prod_{t\in F}S_t\x\prod_{t\in T\sm F}\{1_t\})\cdot U=\prod_{t\in T}G_t$, and
$\card{\prod_{t\in F}S_t}=\max_{t\in F}\bd{G_t}$.
Thus, $\bd{\prod_{t\in T}G_t}\leq\sup_{t\in T}\bd{G_t}$.
For each index $t\in T$, the group $G_t$ is a subgroup of the group $\prod_{t\in T}G_t$. Thus~\cite[Lemma~2.9]{PvK}, $\sup_{t\in T}\bd{G_t}\leq \bd{\prod_{t\in T}G_t}$.
\epf

Recall from the introduction that the compact weight of a topological space is the supremum of the weights of its compact subsets.

\blem\label{lem:LieGen}
Let $X$ be a Tychonoff space of compact weight $\kappa$, and $G$ be an infinite topological group.
\be
\item $\bd{\C(X,G)}\leq \kappa\cdot\wght(G)$.
\item If the group $G$ contains an arc, then $\bd{\C(X,G)}\geq\kappa$.
\ee
In particular, if the group $G$ contains an arc and $\wght(G)\leq \kappa$, then $\bd{\C(X,G)}=\kappa$.
\elem
\bpf
(1)
The group $\C(X,G)$ is topologically isomorphic to a subgroup of the group $\prod_{K\in \K(X)} \C(K,G)$:
The map
\begin{align*}
\psi\colon \C(X,G) &\longrightarrow \prod_{K\in \K(X)} \C(K,G),\\
f &\longmapsto \prod_{K\in \K(X)} f|_K.
\end{align*}
is a topological isomorphism onto its image.
Since the group $\C(X,G)$ is embedded in the group $\prod_{K\in\K(X)}\C(K,G)$, we have
\[
\bd{\C(X,G)}\leq \Bd{\prod_{K\in\K(X)} \C(K,G)}\text{\cite[Lemma~2.9]{PvK}}.
\]
For each compact set $K$ in $X$, we have $\bd{\C(K,G)}\leq\wght(\C(K,G))$~\cite[Proposition~2.11(2)]{PvK}.
By Lemma~\ref{lem:bdProd}, we have
\[
\Bd{\prod_{K\in\K(X)} \C(K,G)}=\sup_{K\in\K(X)}\bd{\C(K,G)}\leq \sup_{K\in\K(X)}\wght(\C(K,G)).
\]
For each compact set $K$ in $X$,
we have $\wght(\C(K,G))\leq \wght(G)\cdot\wght(K)$~\cite[Theorem~3.4.16]{eng}.
Thus,
\[
\sup_{K\in \K(X)}\wght(\C(K,G))\leq \wght(G)\cdot \sup_{K\in\K(X)}\wght(K)=\wght(G)\cdot\kappa.
\]

(2)
Let $A\sub G$ be an arc
containing the identity element $1_G$ of the group $G$.
Take a symmetric open neighborhood $W$ of $1_G$ such that
$A\sm W^2\neq\emptyset$.
Let $K$ be a compact set in $X$.
Let $S\sub \C(X,G)$ be a set of cardinality $\bd{\C(X,G)}$ such that $S\cdot [K,W]=\C(X,G)$.

The family $\cB:=\set{f\inv[W]\cap K}{f\in S}$ is a base for the topological space $K$:
Let $U$ be an open neighborhood in the space $X$ of a point $y\in K$.
Let $a\in A\sm W^2$.
Since the set $A$ is homeomorphic to the unit interval, and $1_G,a\in A$,
there is a function $g\in\C(X,G)$ such that $g[X\sm U]=\{a\}$ and $g(y)=1_G$.
Since $S\cdot [K,W]=\C(X,G)$, there is a function $f\in S$ such that $g\in f\cdot [K,W]$.
Let $\tilde{f}$ be the inverse element of the element $f$ in the group $\C(X,G)$.
We have $\tilde{f}\cdot g\in [K,W]$.
Since $g(y)=1_G$, we have
$\tilde{f}(y)=\tilde{f}(y)\cdot 1_G=\tilde{f}(y)\cdot g(y)=(\tilde{f}\cdot g)(y)$.
Thus $\tilde{f}(y)\in W$, and $y\in \tilde{f}\inv[W]\cap K$.
Assume that there is an element $x\in (\tilde{f}\inv[W]\cap \K)\sm U$.
Then $\tilde{f}(x)\cdot a=\tilde{f}(x)\cdot g(x)=(\tilde{f}\cdot g)(x)$, and thus the element $\tilde{f}(x)\cdot a$ belongs to the set $W$.
Since $\tilde{f}(x)\in W$ and the set $W$ is symmetric, we have $a\in W^2$, a contradiction. Thus, $\tilde{f}\inv[W]\cap K\sub U$.

We conclude that $\bd{\C(X,G)}=\card{S}\geq \card{\cB}\geq \wght(K)$.
\epf

\bprp\label{prp:chiCXG2}
Let $X$ be a Tychonoff space, and $G$ be a topological group containing an arc.
Then $\chi(\C(X,G))=\cof(\K(X))\cdot\chi(G)$.
\eprp

\bpf
By Lemma~\ref{lem:chiCXG}, we have $\chi(\C(X,G))\leq\cof(\K(X))\cdot\chi(G)$.
Let $\cK$ and $\cU$ be directed subsets of the sets $\K(X)$ and $\nbe(G)$, respectively,  such that the set $\sset{[K,U]}{K\in\cK, U\in\cU}$ is cofinal in the set $\nbe(\C(X,G))$, of cardinality $\chi(\C(X,G))$.
By Proposition~\ref{prp:chiCXG}, we have $\cof(\K(X))\leq\card{\cK}$, and $\chi(G)\leq\card{\cU}$.
Thus, $\cof(\K(X))\cdot \chi(G)\leq\chi(\C(X,G))$.
\epf

\bpf[Proof of Theorem~\ref{thm:LieGen}]
The inequality
\[
\bd{\C(X,G)}\leq \ds(\C(X,G))\leq\bd{\C(X,G)}\cdot \chi(\C(X,G))\text{~\cite[Proposition~2.11(2)]{PvK}}
\]
implies that $\bd{\C(X,G)}\cdot \chi(\C(X,G))=\ds(\C(X,G))\cdot \chi(\C(X,G))$.
By $\wght(\C(X,G))=\ds(\C(X,G))\cdot \chi(\C(X,G))$, and Proposition~\ref{prp:chiCXG2}, we have
\[ \wght(\C(X,G))=\bd{\C(X,G)}\cdot\cof(\K(X))\cdot\chi(G).
\]
By Lemma~\ref{lem:LieGen}(1), since the compact weight of a metric space is countable, we have $\bd{\C(X,G)}\leq \wght(G)$. Thus,
\[ \cof(\K(X))\leq\wght(\C(X,G))=\cof(\K(X))\cdot\wght(G).
\]
Apply Theorem~\ref{thm:hi5}.
\epf

Theorem~\ref{thm:LieGen} implies the following result.
\bcor\label{cor:wghtCXG}
Let $X$ be a complete metric nondiscrete noncompact space, and $G$ be a second countable
topological group containing an arc.
\be
\item If the space $X$ is locally compact, then $\wght(\C(X,G))=\ds(X)$.
\item If the space $X$ is not locally compact, then
\[
\wght(\C(X,G))= \ds(X)\cdot\fd\cdot\sup_{x\in X}\Cof{[\ld(x)]^\alephes}.\qed
\]
\ee
\ecor

\brem
For a metric space $X$, let $A(X)$ be the free abelian topological group generated by the space $X$.
Let $\widehat{A(X)}$ be the dual group to the group $A(X)$.
Let $\mathbb{T}$ be the one-dimensional torus.
By results of Pestov~\cite[Proposition~7]{pestov}, and Galindo--Herna\'ndez~\cite[Theorem~2.1]{GH}, the
topological group $\widehat{A(X)}$ is topologically isomorphic to $\C(X,\mathbb{T})$.
Thus, Theorem~\ref{thm:exactGbase} and Corollary~\ref{cor:wghtCXG} apply to the group $\widehat{A(X)}$.
\erem

\end{document}